\documentclass[11pt,a4paper]{article}

\usepackage{amsfonts}
\usepackage{graphics}
\usepackage{amssymb}
\usepackage{amsmath,amsthm,amsfonts,amssymb,amscd}
\usepackage[dvips]{graphicx}

\newtheorem{theorem}{Theorem}[section]
\newtheorem{proposition}[theorem]{Proposition}
\newtheorem{corollary}[theorem]{Corollary}
\newtheorem{remark}[theorem]{Remark}
\newtheorem{definition}[theorem]{Definition}
\newtheorem{lemma}[theorem]{Lemma}

\newcommand{\eop}{\hfill $\Box$}

\newcommand{\Prob}{\mathbf P}
\newcommand{\R}{\mathbf R}
\newcommand{\N}{\mathbf N}

\newcommand{\Z}{\mathbf Z}
\newcommand{\E}{\mathbf E}

\newcommand{\fL}{\mathcal L}

\begin{document}

\begin{center}

\vspace{1cm}

 {\Large {\bf Degenerate semigroups and stochastic flows of mappings  in
foliated manifolds\\[3mm]
}}

\end{center}

\vspace{0.3cm}

\begin{center}
{\large {Paulo Henrique P. da Costa}
\footnote{e-mail: phcosta@ime.unicamp.br.
} \ \ \ \ \
\ \ \ \ \ \ \ \   Paulo R.
Ruffino}\footnote{Corresponding author, e-mail:
ruffino@ime.unicamp.br.
}

\vspace{0.2cm}

\textit{Departamento de Matem\'{a}tica, Universidade Estadual de Campinas, \\
13.083-859- Campinas - SP, Brazil.}

\end{center}

\begin{abstract} Let $(M, \mathcal{F})$ be a compact Riemannian foliated 
manifold.
We consider a family of compatible Feller
semigroups in $C(M^n)$  associated to laws of the $n$-point 
motion.
Under some assumptions (Le Jan and Raimond, \cite{Le Jan-Raimond}) 
there
exists a stochastic flow of measurable mappings in $M$. We study the
degeneracy of these semigroups such that the  
flow of mappings is
foliated, i.e. each trajectory lays in a single leaf of the foliation a.s,
hence creating a geometrical obstruction for coalescence of trajectories in
different
leaves. As an application, an
averaging principle is proved for a first order perturbation transversal to the
leaves. Estimates for the rate of
convergence are calculated.
\end{abstract}


\vspace{0.3cm}
\noindent {\bf MSC2010 subject classification:} 60H10, 60H30, 57C12.

\section{Introduction}

The understanding of the geometry and dynamics in a foliated 
space plays a quite important role in the intersection of many areas. Among 
many others contemporary good references, we mention few of them which are 
closer to the approach we are interested in this paper, e.g.  Tondeur 
\cite{Tondeur}, Candel and Conlon \cite{candel}, Plante\cite{plante},  Walczak 
\cite{Walcak}. On what regards stochastic systems, their interlace with 
foliations received a boost with the paper by L. Garnnet \cite{Garnett}, where 
she has constructed harmonic measures introducing foliated Brownian motions on 
the leaves. Since them, many works intertwinning stochastics, semigroups and 
foliations  have appeared, see e.g. Kaimanovich \cite{kaimainovich}, Candel 
\cite{Candel_Adv_Math}, \cite{Catuogno-Ledesma-Ruffino} and references therein.

\smallskip

In this article we are interested on constructing and studying properties of a 
certain class of stochastic flow adapted to a foliation in the sense that each 
trajectory lays in a single leaf of the the foliation. Classically, for 
stochastic (Stratonovich) 
differential equation in a differentiable manifold with smooth vectors fields, 
the existence of a flow of (local) diffeomorphisms has been well studied and 
has 
been applied in a vast range of topics in the literature. We mention mainly 
Kunita \cite{Kunita 1} \cite{Kunita 2} and  \cite{L Arnold} 
which contain in the 
references many other authors that contributed in the construction of these 
flows, each of them has considered different relevant aspects and properties. 
Extending this classical flow with regularity in the space and continuity of 
trajectories with respect to time $t\geq 0$, many constructions have 
appeared, see e.g.   a survey on some of these 
nonclassical stochastic flow in Tsirelson \cite[In particular, 
Chap.7]{Tsirelson}, also \cite{Baxendale}, \cite{Bertoin-Le Gall}, 
\cite{Chen_Xiang} among others.  In this article we are particularly interested 
on studying degenerate dynamical properties of the
stochastic flow of mappings as constructed by Le Jan and Raimond
\cite{Le Jan-Raimond}, where coalescence also can happen.

 \smallskip

 The Le Jan-Raimond stochastic flow, which here is going to called by the 
acronym  LJR-flow, is 
constructed over a (locally) compact metric space $M$, based
on a family of
Feller compatible semigroups of operators in $M^n$ for integers $n\geq 1$ . 
That is, for each $n \geq 1$, there exists a semigroup which determines the law 
of a Markovian process of the $n$-point motion. Given compatibility and 
diagonal preserving 
conditions on these semigroups, there exists a stochastic flow of measurable 
mappings in $M$ which generate this family of semigroups. These concepts and 
results are going to be precisely stated and recalled 
below.

\smallskip

In this article we consider $M$ a Riemannian compact manifold.
 Motivated by the 
concept of invariant submanifolds for flows, we look initially for conditions 
on the family of Feller semigroups which 
guarantee that a certain submanifold $N$ is invariant by the corresponding 
LJR-flow. In this case, we call the semigroups $N$-degenerate, 
cf. Definition \ref{Def: sg N degenerado}. With appropriate hypothesis, Theorem 
\ref{Thm: 
fluxo_foleado} guarantees that trajectories starting in $N$ lays in $N$ 
 for all time 
$t\geq 0$ a.s.. Suppose now that the manifold $M$ has a regular foliation 
$\mathcal{F}$. 
Following the ideas ideas for a submanifold, an 
appropriate degeneracy on the 
semigroups with respect to the leaves  of 
the foliation  (Definition \ref{Def: sg folheado}) will
imply that each trajectory lays in a single leaf of the foliation, as stated 
in  Theorem \ref{Thm: fluxo_foleado}. This foliated behaviour of the flow, in 
particular, introduces a
geometrical obstruction to coalescence of points in different leaves.

\smallskip

The class of examples on which this foliated phenomenon happens includes 
stochastic dynamical
systems generated by classical
Stratonovich differential equations whose vector fields belong to
the tangent space of the leaves. Degenerate stochastic systems in
differentiable manifolds also illustrates
this context if the dimension of the H\"{o}rmander Lie algebra space has
constant dimension, hence, in this case the foliation is determined a posteriori
by the vector fields of the original SDE.

\medskip

The article is organized in the following way. In the next  paragraphs we 
recall the
concepts and results on the construction of the
LJR-flow, as in
\cite{Le Jan-Raimond}, which we are going to use in the next sections. In
Section 2 we introduce the definitions of
$N$-degenerate and foliated semigroups. In Section 3 we proof the main result 
on foliated flow in $M$ which guarantees the equivalence between the foliated 
semigroups and LJR foliated flows. The reader will notice that for 
noncompact $M$, the same result holds if one assume some further, rather 
natural  conditions (see also \cite{Le Jan-Raimond}). Coalescent foliated flow 
is discussed at the 
end of this section. Finally in Section 4, as an
application of this approach, assuming that the leaves are compact, we 
investigate the effective behaviour  of a small 
transversal perturbation of order $\epsilon$ in a foliated LJR-flow which 
destroys the foliated 
structure. An averaging principle is proved for a perturbation given by a 
vector field
$\epsilon K$, with $K$ transversal to the leaves. We assume hypothesis on 
convergence of average functions along the perturbed 
process (Hypothesis H1 and H2) which are natural in most stochastic systems, 
see \cite{ivan, Hoegele-Ruffino, Li}. Essentially, for small  
$\epsilon>0$, the transversal behaviour, with time rescaled by 
$\frac{t}{\epsilon}$, is approximate by an ODE in this space whose vector field 
is given by the ergodic average of the transversal component of $K$ on each 
leaf. We find estimates on the rate of
convergence when $\epsilon$ goes to $0$.

\subsection{Le Jan-Raimond stochastic flow of measurable mappings}

We summarize below the main definitions and results on LJR-flows cf.
\cite[section 1]{Le Jan-Raimond}. We apply this theory to our specific context
here, where the state
space $M$ is a Riemmanian manifold. In the next sections $M$
will additionally  be
endowed with a regular foliation. We shall consider Markovian semigroups acting
on $B_b(M^n)$, the space of bounded measurable functions in $M^n$ and
Feller continuous semigroups in $C_0(M^n)$, the space of continuous functions in
$M^n$ which goes to zero at infinity, see e.g. among others Revuz and Yor 
\cite{Revuz-Yor}.

The results in the next sections are based on the fact that a consistent
system of
$n$-point Feller semigroup (which can be viewed as the law of the motion of $n$
indivisible points thrown into a fluid in $M$) and a preserving condition on
the diagonal determine uniquely in law the flow of measurable mappings in
$M$. Next paragraphs make this statement more precise, essentially we have that
the next three definitions are equivalent, cf. Theorem \ref{Thm_fluxo_LJR}.

\begin{definition} \label{Def: familia CDP}
\emph{Let $ ( P_{t}^{(n)}, n \in \N )$ be a family of Feller
(Markovian) semigroups
in $C_0(M^n)$
 (in $B_b(M^n)$, respectively).}
    \begin{enumerate}
     \item \emph{The family $( P_{t}^{(n)} )$ is called
\textit{compatible} if
      for all  $f \in C_0(M^n)$ which can be written in terms of fewer
variables, i.e. there exists
      an $\bar{f}\in C_0(M^k)$, $k< n$,  with  
$f(x_1,\cdots,x_n) \, = \,
      \bar{f}(x_{i_1},\cdots,x_{i_k})$, 
      we have that
      \[
       P_{t}^{(n)}\, f(x_1,\cdots,x_n)= P_{t}^{(k)} \,
\bar{f}(x_{i_1},\cdots,x_{i
       _k}).
      \]}
     \item \emph{The family $( P_{t}^{(n)}) $ is \textit{diagonal
preserving} in the sense that for all $f \in C_0(M)$
            \[
                P_{t}^{(2)} f^{\otimes 2}(x,x) \ = \ P_{t}\  f^2(x)
            \]where here and along the article, $P_t$ stands for $P_t^{(1)}$.}
    \end{enumerate}
\end{definition}

We shall abbreviate and say that a family of Feller semigroups $(
P_{t}^{(n)}) $ is CDP if it is compatible and
diagonal preserving as defined above.
For a fixed $n\in \N$, the Markov process associated to $P_t^{(n)}$ starting at
a certain initial condition $(x_1, \ldots, x_n)$ is called the $n$-point motion
of this family of semigroups, it is defined on the set of c\`adl\`ag paths on
$M^n$. See more on $n$-point motion also in Kunita \cite{Kunita 2}.
\bigskip

Let  $(F,\mathcal{E})$ be the space of
measurable mappings on $M$ endowed with the $\sigma$-algebra generated by
the
application map $\varphi\longmapsto \varphi(x)$ for every $x\in M$.
More precisely, $\mathcal{E} = \sigma \left\{ \{ \varphi \in
F: \varphi(x) \in A \}, \mbox{ for all } A \in  \mbox{Borel}(M),
\mbox{ and } x \in M \right\} $.

\begin{definition} \label{semigroup_convolution}
\emph{ Consider a convolution semigroup $\{Q_t: t \geq 0\}$ of regular
probability measure on
$(F,\mathcal{E})$, i.e $Q_{s+t} =
    Q_s \ast Q_t$ for all $0 \leq s, t$. This family is called a Feller
convolution semigroup  if for all $f \in C_0(M)$ we have that:
    \begin{enumerate}
        \item[a)] \ \ $\displaystyle \lim_{t \to 0} \sup_{x
\in M} \int
        (f\circ \varphi (x) - f(x))^2\,Q_t(d\,\varphi) = 0;
        $
        \item[b)] For all $x\in M$ and $t\geq 0$
\[
\lim_{y\rightarrow x} \int (f\circ
        \varphi (y) - f\circ \varphi(x))^2\,Q_t(d\,\varphi) = 0;
\]
and
\[
 \lim_{y\rightarrow \infty} \int \left(f \circ
        \varphi (y) \right)^2\, Q_t (d\, \varphi)=0.
\]
    \end{enumerate}
}
\end{definition}

The third and last equivalent definition in this context is the following. Let
$(\Omega, \mathcal{F}, \Prob)$ be a probability space.

\begin{definition} \label{Def: stochastic_flow}\emph{
A family of $(F, \mathcal{E})$-valued random
variables $\varphi = (\varphi_{s,t}, 0\leq s \leq t)$  is called a
\textit{measurable stochastic flow of mappings} if $(x, \omega) \mapsto
\varphi_{s,t}(x)$ is measurable, stationary, with independent increments,
satisfies the cocycle property, i.e. for all $0\leq s\leq u \leq t$ and $x\in
M$, $\Prob$-a.s. $\varphi_{s,t}= \varphi_{u,t}\circ \varphi_{s,u} (x)$ and for
every $f \in C_0(M)$ we have that:
    \begin{enumerate}
        \item[a)] \ \ $\displaystyle \lim_{(u,v) \rightarrow (s,t)} \ \
\sup_{x
\in M} \E
         \left[ \left( f\circ \varphi_{u,v} (x) - f\circ
\varphi_{s,t} (x)\right)^2 \right]= 0;
        $
        \item[b)] For all $x\in M$ and $t\geq 0$
\[
\lim_{y\rightarrow x} \E
\left[ \left( f\circ  \varphi_{0,t} (y) - f\circ \varphi_{0,t}(x)\right)^2 
\right] = 0;
\]
and
\[
 \lim_{y\rightarrow \infty}
 \E \left[ \left( f \circ \varphi_{0,t} (y) \right)^2 \right] =0.
\]
    \end{enumerate}
A family of $(F, \mathcal{E})$-valued random
variables $\varphi = (\varphi_{s,t}, 0\leq s \leq t)$  is called a
\textit{stochastic flow of mappings}
if there exists a measurable stochastic flow of mappings $\varphi'=
(\varphi'_{s,t}, 0\leq s \leq t)$ such that for all $x\in M$ and $0\leq s
\leq t$, we have  that  $\varphi_{s,t}(x)= \varphi'_{s,t}(x)$ $\Prob$-a.s.  .
}
\end{definition}

Given a stochastic flow of mappings $\varphi_{s,t}$,  the
law of $\varphi_{0,t}$ determines a semigroup of convolution $Q_t$; which in
turn  determines a family of Feller CDP semigroups given by, for $f\in C_0(M^n)$
and $x\in M^n$,
\begin{equation}  \label{CDP_convolution}
 P_t^{(n)} f(x) = \int_F f\circ\varphi^{\otimes
    n}(x)\ Q_t (d\varphi),
\end{equation}
see Propositions 1.2 and 1.3 in \cite[p.1252-1254]{Le Jan-Raimond}. In fact the
converse also hold, hence Definitions \ref{Def: familia CDP},
\ref{semigroup_convolution} and \ref{Def: stochastic_flow} are associated 
one to each other, as stated
in
the next theorem.

\bigskip

Denote by $(\Omega, \mathcal{A})$ the space $(\Pi_{s\leq t}
F,\bigotimes_{s\leq t} \mathcal{E})$ and by
$\varphi_{s,t}$ the canonical stochastic process $\omega \mapsto
\omega(s,t)$.

\begin{theorem}[Le Jan-Raimond] \label{Thm_fluxo_LJR} A family of Feller CDP 
semigroups
$( P_{t}^{(n)}, n\in \N)$ determines a unique Feller convolution semigroup
of probability measures $(Q_t)_{t\geq 0}$ on $(F, \mathcal{E})$ which satisfies
Equation (\ref{CDP_convolution}). A Feller convolution semigroup
$(Q_t)_{t\geq 0}$ in turn determines a unique shift invariant probability
measure $\Prob_Q$ on
$(\Omega, \mathcal{A})$ such that the canonical stochastic
process $(\varphi_{s,t}: s \leq t)$ is a
stochastic flow of mappings with law $Q_{t-s}$.
\end{theorem}


A stochastic flow of measurable mappings on $M$, $(\varphi_{s,t}: s\leq
t)$, is
called a coalescing flow if for some point $(x, y) \in M^2$, $T_{x,y} =
\inf \{t\geq 0: \varphi_{0,t} (x) \, = \, \varphi_{0,t} (y)\}$ is
finite with positive probability and for any $t \geq
T_{x,y}$
we have that $\varphi_{0,t} (x) \, = \, \varphi_{0,t} (y)$. In other words, a
stochastic flow is coalescing if there exists $x\neq y$ such that their
trajectories stick together after a finite time with positive probability.


\section{Definitions and Preliminary results}

Consider a compact Riemannian manifold $M$. We shall use the extended 
definition of a submanifold of $M$, as in Candel and Conlon \cite{candel}, in 
the 
following sense: we say that $N\subset M$ is a submanifold when 
$N$ is the image of a connected manifold $W$ by an injective immersion  $i:W 
\rightarrow 
M$, that is $N=i(W)$. Whence $N$ is endowed with two 
topologies: the intrinsic topology $\tau_{\mathrm{int}}$ of $W$ and the 
induced topology $\tau_{\mathrm{ind}}$ as a subset of $M$. Obviously 
$\tau_{\mathrm{ind}} 
\subseteq \tau_{\mathrm{int}}$, and equality holds if $i$ is an embedding. Most 
of 
the time we are interested in the case of $N$ being a complete submanifold with 
respect to the induced Riemannian metric. In this case 
we have that 
\[
 C_0(N, \tau_{\mathrm{int}})= \left\{ f|_N; f\in C(M) \right\}
\]
if and only if $N$ is compact with respect to $\tau_{\mathrm{int}}$, i.e.,  $W$ 
is 
compact. Hence, in this case, we have also that $\tau_{\mathrm{int}}= 
\tau_{\mathrm{ind}}$ and 
$C_0(N, \tau_{\mathrm{int}}) = C(N, \tau_{\mathrm{int}})$. Note that if $N$ is 
dense in
$M$ then $C_0(N, \tau_{\mathrm{int}}) \cap C(N, \tau_{\mathrm{ind}})= \{0\}$.

\bigskip

The geometrical idea of invariant manifolds for certain dynamics, say 
e.g. stable submanifold,  motivates the following definition.

\begin{definition}\label{Def: sg N degenerado}
\emph{Let $N$ be  a submanifold of $M$, a Feller
semigroup
$P_{t}$ is $N$-degenerate if there exists a
Markovian semigroup $P_{t}^{N}$ such that
for all $f \in C(M)$ we have that:
\begin{equation} \label{Eq: da def de sg N-degenerado}
		P_{t} f (x) \ = \ P_{t}^{N}
        f\big|_{N}(x),
\end{equation}
for all $x \in N$.}
\end{definition}
Since the set of test functions for the definition above reduces to a subset of 
$C(N, 
\tau_{\mathrm{ind}})$, the Markovian semigroup $P_t^N$ can 
be non unique in $C_0(N, \tau_{\mathrm{int}})$. If $N$ is compact then 
$P_t^N$ is unique and Feller. 

\paragraph{Foliations.} A regular foliation $\mathcal{F}$ in $M$ is a partition
of $M$ into equivalent classes of complete submanifolds of the same dimension;
it corresponds to assign a regular integrable differentiable $d$-dimensional
distribution in $M$. Each of these submanifolds are called the leaves of the
foliation $\mathcal{F}$. More precisely: an  
 $(n+d)$-dimensional smooth manifold $M$ is foliated when there exists
an
atlas on $M$  such that for any pair $(U,\psi)$ and $(V,\phi)$ of
coordinate maps we have that:
\begin{enumerate}
\item  $\psi(U) \,=\,U_1\,\times\,    U_2$,
    where $U_1 \subseteq \R^n$ and $U_2 \subseteq \R^{d}$ are open sets;
\item If   $U \cap V \neq \emptyset$ then the map 
$\phi\circ\psi^{-1}:\psi(U \cup V)  \rightarrow \phi(U \cup V)$ has the
following form
$\phi\circ\psi^{-1}(x,y)    \, = \, (h_1(x,y),h_2(y))$.
\end{enumerate}
Given a point $x\in M$, the leaf passing through $x$ is denoted by
$\mathcal{L}_x$. For further properties and details see e.g 
\cite{camacho},  \cite{candel}, \cite{Tondeur},  \cite{Walcak}.

\begin{definition} \label{Def: sg folheado}
\emph{Let $(M, \mathcal{F})$ be a foliated space. A Feller semigroup
$P_{t}$ in $M$ is {\it 
foliated} or 
$\mathcal{F}$-{\it foliated} if  $P_{t}$ is $\fL$-degenerate for every $\fL 
\in
\mathcal{F}$, i.e. there exists a family of Markovian 
semigroups 
$\left(
P_t^{\mathcal{L}}\right)_{\mathcal{L}\in \mathcal{F}}\ $ such that for all $f
\in C(M)$
\begin{equation} \label{Eq: da definicao sg folheado}
		P_{t} f (x) \ = \ P_{t}^{\mathcal{L}_ x}
        f\big|_{\mathcal{L}_ x}(x),
\end{equation}
for all $x \in M$.}
\end{definition}


As before, if the leaves are compact then for each leaf $\fL$, the semigroup 
$P_t^{\mathcal{L}}$ which satisfy Equation (\ref{Eq: da definicao sg folheado}) 
is Feller and unique. In this case, there exists also an intrinsic way
to verify  whether a certain Feller semigroup $P_t$ in $M$ is foliated 
(degenerate) or
not.
For this purpose, given $f \in C(M)$ and a submanifold $N\subset M$,
consider the set 
    \[
        I_{f,N} \ := \ \left\{g \, \in C(M):
        g{\big|}_{N} \ = \ f{\big|}_{N}\right\}.
    \]
Hence, by definition, if $P_t$ is an $N$-degenerate semigroup then for
any $g \in I_{f,N}$ we have $
        P_{t} \, f(x) \ = \ P_{t} \, g(x) $. Reciprocally:

\begin{proposition} Consider a Feller semigroup $P_t^M$ on $C (M)$.

\begin{enumerate}
    \item Assume that the submanifold  $N$ is compact. The semigroup  
$P_{t}$ is $N$-degenerate if and only if $P_{t}
\, f(x) \
= \
        P_{t}\, g(x)$ for all $g \in I_{f,N}$.
    \item Assume that the leaves of $\mathcal{F}$ are compact. The semigroup 
$P_{t}$ is foliated if and only if
$P_{t} \, f(x) \
= \
        P_{t}\, g(x)$ for all $g \in I_{f,\mathcal{L}}$ and all
$\mathcal{L}\in \mathcal{F}$.
\end{enumerate}
\end{proposition}

\proof We prove the converse of item (1). We only have to show that 
there exists a
Feller semigroup acting on $C(N)$ which satisfies
Equation (\ref{Eq: da def de sg N-degenerado}). In fact, given $g \, \in C(N)$, 
take an  extension  $f \in C(M)$ and 
define $P_t^N g (x)= P_t f (x)$ for all $x\in N$. Hypothesis guarantees that 
$P_t^N$ is well defined. Item (2) follows trivially now. 

\eop

%
%
%
%

Contrasting with last Proposition, if $N$
is dense in $M$ then each set $I_{f,N}$ has a
single element.

\section{Foliated flow}

At this point we return to the family $( P_{t}^{(n)}, n\in \N) $  of 
Feller CDP
semigroups as defined in Section 1.1. Initially note that if 
$\mathcal{F}$ is a
foliation of $M$, then $\bar{\mathcal{L}} =\mathcal{L}_1\times \ldots \times 
\mathcal{L}_n$ with
$\mathcal{L}_j \in \mathcal{F}$, $j=1, \ldots,  n$, are leaves of a foliation 
in 
$M^n$ which we are going to denote by $\mathcal{F}^n$. The main results in 
this
section show that
for a family of Feller CDP semigroups $( P_t^{(n)}, n\in \N)$, if
$P_t^{(i)}$ is $\mathcal{F}^i$-foliated for a certain $i\in \N$, then
$P_t^{(n)}$ is $\mathcal{F}^n$-foliated for all $n\in \N$, moreover the
associated LJR-flow $\varphi_{s,t}$
is foliated in the following sense: 

\begin{definition}\label{Def: foliated stoch flow}
\emph{A stochastic flow of mappings $\varphi_{s, t}$ in $M$
is called an $\mathcal{F}$-foliated stochastic flow if  for all $0\leq s,t$ we
have that $\varphi_{s, t} (x) \in \mathcal{L}_x$ a.s..}
\end{definition}

%


The result below states that if the lowest level Feller semigroup $P_t^{(1)}$ is
$\mathcal{F}$-foliated then the corresponding LJR-flow  of mappings 
$\varphi_{s,t}$
is $\mathcal{F}$-foliated.

\begin{theorem}[Foliated flow] \label{Thm: fluxo_foleado} Let 
$( P_{t}^{(n)}: n \in
\N ) $ be a family of 
Feller CDP semigroups in $M$.
\begin{enumerate}
\item\emph{\textbf{(Invariant submanifold)}} Let $N$ be a submanifold
of $M$. If $P_{t}^{(1)}$ is $N$-degenerate then the LJR stochastic flow of
mappings preserves $N$, i.e.  for  $0\leq s \leq t$ and  $x\in N$
we have that $\varphi_{s,t}(x)\in N$,  $\Prob$-a.s..
\item\emph{\textbf{(Foliated Flow)}} If $P_{t}^{(1)}$ is $\mathcal{F}$-foliated
then the corresponding LJR
    stochastic flow of mappings $\varphi_{s,t}$ is an $\mathcal{F}$-foliated
flow.
\end{enumerate}
\end{theorem}

\proof For item (1), consider an increasing sequence  of compact
sets with respect to the intrinsic topology $\tau_{\mathrm{int}}$ in $N$ such 
that $N =
\cup_n K_n$. The idea of the proof is to control the probabilities with which
the process exits the sets $K_n$. For
fixed $0 \leq t$ and $x
\in N $ consider the measurable set
\[
 B \ = \ \left\{\omega \, \in \Omega:
                \varphi_{0,t}(x,\omega) \, \in \ M \setminus N
\right\}.
\]
We prove that $\Prob (B)=0$ writing $ B =  \cap_n
D_n$ with
\[
D_n \ = \ \left\{\omega \, \in \Omega:
                \varphi_{0,t}(x,\omega)  \in M \setminus K_n \right\}
\]
and proving that $\Prob (D_n)$ goes to zero when $n$ tends to infinity. The 
semigroup  $P_{t}^{N}$ has an associated transition probability measure
$\mu_{P_{t}^{N}}(x, dy)$ with support in $(N, \tau_{\mathrm{int}})$, see e.g. 
Revuz and Yor
\cite[Chap. III.2]{Revuz-Yor}, which in general does not coincide with the 
support of $\mu_{P_{t}^{(1)}}(x, dy)$ in $M$. 
The key point here is to link these two 
probability 
measures using continuous functions in $C(M)$, as demanded in Definition 
\ref{Def: sg N degenerado}.

Consider closed sets
$F_{j,n} \subset M \setminus K_n$ which are increasing in $j$ and such that $M -
K_n =
\cup_j F_{j,n}$. For each pair $(j,n)$ take a continuous function $f_{j,n
} \in C(M, [0,1])$ such that
$f_{j, n}(x) =  1$ for all
$x \in F_{j,n}$ with support in $ M \setminus
 K_n$.
 Therefore, for a fixed $x \in K_n$, by formula (\ref{CDP_convolution}) and 
the fact that $P_t$ is $N$-degenerate we have that 
    \begin{eqnarray*}
            \E\left[ f_{j,n}\circ \varphi_{0,t}(x) \right] & =&
P_{t}^Nf_{j,n} \big|_N (x) \\
&& \\
&=&
\int_N f_{j,n} \big|_N (y) \ \mu_{P_{t}^N}(x, dy).
    \end{eqnarray*}
Hence, using the fact that $f_{j,n}$ converges pointwise to the characteristic
function $1_{M\setminus K_n}$ when $j$ goes to infinity we 
conclude that
\begin{eqnarray*}
\Prob (D_n) & = &
    		\lim_{j \to \infty}\E\left(f_{j,n} \circ
            \varphi_{0,t}(x)\right)\\
 && \\
 &=&  \int_N \lim_{j \to \infty} f_{j,n}
            \big|_N (y) \ \mu_{P_{t}^N}(x, dy) \\
 &   & \\
 &= & \mu_{P_{t}^{N}}(x, N\setminus K_n).
    \end{eqnarray*}
Last term  goes to zero when $n$ goes to infinity since
$\mu_{P_{t}^{N}}(x,
dy)$ is a Radon measure.

Item (2) of the statement follows directly from item (1) applied in each leave
of the foliation.

\eop

%
%

If the submanifold $N$ in the item (1) of Theorem \ref{Thm:
fluxo_foleado} above is compact, a proof purely analytical functional on $C(M)$ 
shows nuances of the technique: 
there exists a
countable sequence of closed sets $F_n$ such that
$\cup_n  F_n = M \setminus N$. For each $n \in \N$, consider a corresponding 
continuous function $f_n\, \in C(M, [0,1])$ with support in $M \setminus
N$ and such that  $f_n(x)  =  1$ for all $x$ in $F_n$. For fixed $0 \leq
t$ and  $x \in N $ we show that the measurable set
\[
 B \ = \ \left\{\omega \, \in \Omega:
                \varphi_{0,t}(x,\omega) \, \in \ M \setminus N
\right\}
\]
has probability zero. Now we write $B = \cup_n
B_n$ where
\[
B_n \ = \ \left\{\omega \, \in \Omega:
                \varphi_{0,t}(x,\omega) \, \in F_n \right\}.
\]
We prove that $\Prob (B_n)= 0 $ for every $n \in \N$. In fact, by 
Chebyshev inequality, Theorem
$\ref{Thm_fluxo_LJR}$  and the definition of foliated semigroup we calculate:
\begin{eqnarray*}
 \Prob (B_n) &\leq&  \E(f_n\circ\varphi_{0,t}(x))\\
    & = &  P_{t}^{(1)}f_n (x) \\
    & = &  P_t^{N} f_n  \bigg|_N = 0.
 \end{eqnarray*}

 \bigskip

 Next result exploits the fact that, although in Theorem \ref{Thm: 
fluxo_foleado}, for simplicity, we have assumed that the first element of the 
family, i.e. with $n=1$, $P^{(1)}_t$ is foliated, the same result holds if one 
assumes, instead, that  $P^{(k)}_t$ is $\mathcal{F}^k$-foliated for a certain 
$k \geq 1$. 


\begin{corollary} \label{Cor: equivalencias de folheacoes}
Let 
$( P_{t}^{(n)}: n \in
\N )$ be a
Feller CDP-family
of semigroups in $C(M)$.
\begin{enumerate}
\item\emph{\textbf{(Invariant submanifold)}} If for a positive integer $k\geq 
1$ 
the semigroup  
$P_t^{(k)}$ is $N^k$-degenerate, then all members of the family 
$(P_{t}^{(n)}: n \in \N )$ are $N^n$-degenerate.

\item\emph{\textbf{(Foliated Flow)}} If for a positive integer $k\geq 1$ the 
semigroup  
$P_t^{(k)}$ is $\mathcal{F}^k$-foliated, then all members of the family 
$(P_{t}^{(n)}: n \in \N )$ are foliated in the  
corresponding $\mathcal{F}^n$-foliation of $M^n$.

\end{enumerate}

\end{corollary}

 \proof We prove Item (2). 
Let $\varphi_{s,t}$ be the LJR-flow associated to 
the CDP-family of
Feller semigroups $( P_{t}^{(n)}: n \in \N )$. Assuming that $P_t^{(k)}$ is 
$\mathcal{F}^k$ foliated, fix a leaf $\bar{\mathcal{L}}= \mathcal{L}_1 \times 
\ldots \times \mathcal{L}_k \in \mathcal{F}^k$, and a point 
$x=(x_1, \ldots, x_k) \in
\bar{\mathcal{L}}$. Let $K_n$ be a sequence of compact sets with respect to the 
intrinsic topology in $\bar{\mathcal{L}}$ such that 
$\bar{\mathcal{L}} = \cup_n K_n$. Controlling the probability with which the 
process $\varphi_{0,t}(x)=(\varphi_{0,t}(x_1), \ldots , \varphi_{0,t}(x_k))$ 
exits the sets $K_n$, following the same argument as in the proof of Theorem 
\ref{Thm: fluxo_foleado}, we have that for all $t\geq 0$, $\varphi_{0,t}(x) \in 
\bar{\mathcal{L}}$ a.s.. This implies that the leaves 
in $\mathcal{F}$ are invariant by the flow $\varphi_{s,t}$ a.s.. 

Now, for any 
$n\geq 1$, given a leaf (using the same notation) $\bar{\mathcal{L}} \in 
\mathcal{F}^n$,  a point 
$x=(x_1, \ldots, x_n) \in
\bar{\mathcal{L}}$ and $g \in B_b (\bar{\mathcal{L}})$, define the Markovian
semigroup $P_t^{\bar{\mathcal{L}}}$ in the leaf $\bar{\mathcal{L}}$  by the 
formula
\begin{equation*}
 P_t^{\bar{\mathcal{L}}} g  (x) = \E \left[ 
g  \left( \varphi_{0,t}(x_1), 
\ldots,
\varphi_{0,t}(x_n) \right) \right].
\end{equation*}
For a function $f \in C(M^n)$, we have obviously that 
\begin{equation*}
 P_t^{(n)} f  = \E \left[ 
f |_{\bar{\mathcal{L}}} \left( \varphi_{0,t}(x_1), 
\ldots,
\varphi_{0,t}(x_n) \right) \right] = P_t^{\bar{\mathcal{L}}} f 
|_{\bar{\mathcal{L}}}.
\end{equation*}
Hence  $P_t^{(n)}$ satisfies Definition \ref{Def: sg folheado} of an 
$\mathcal{F}^n$-foliated semigroup.
Item (1) follows directly by the same argument.

\eop

Next corollary establishes sufficient conditions for the semigroups in the 
leaves to be Feller.

\begin{corollary}[Feller in the leaves] \label{Cor: sg Feller na folha}
If the foliated flow of measurable mappings established by 
Theorem \ref{Thm: fluxo_foleado} is such that   $\varphi_{s,t}|_{\mathcal{L}}$ 
satisfies Conditions  
(a) and (b) of 
Definition \ref{Def: stochastic_flow} for all $f\in C_0 (\mathcal{L})$ then 
the semigroups in the leaves $P^{\bar{ \mathcal{L}}}_t$ for $\bar{ 
\mathcal{L}} \in \mathcal{F}^n$ are also Feller for all $n\geq 1$.
(Analogous for the $N$-invariant flows).
\end{corollary}

%
%
%
%


\proof  
By Equation 
$\ref{CDP_convolution}$ and Theorem \ref{Thm_fluxo_LJR}, (alternatively 
\cite[Prop. 1.3]{Le Jan-Raimond}) applied in each leaf $\bar{\mathcal{L}}$, it 
follows that the semigroups $P^{\bar{\mathcal{L}}}_ t$ are Feller. 

\eop

\noindent {\bf Remark:} Compactness of the leaves or continuity of   
$\varphi_{0,t} |_{\mathcal{L}}(x)$ in  $x$ and $t$ with respect to the 
intrinsic 
topology of $\mathcal{L}$ imply Conditions (a) 
and (b) in the hypothesis of Corollary 
\ref{Cor: sg Feller na folha} above.

\paragraph{Example: Foliated semigroups in dense leaves.}

We consider the canonical example of a flat 
2-torus $T = \R^2/\Z^2$ and a unitary vector in the plane $v=(v_1, 
v_2)\in \R^2$. We take a foliation in $T$ such that at each point 
$x= (a,b) \in T$, the leaf passing thorough $x$ is given by the 
winding of a line passing thorough it in direction of $v$:
\[
 \mathcal{L}_{(x,y)}= \{ (a+\lambda v_1, b+ \lambda v_2) \ \ (\bmod \ \Z^2), 
\mbox{ for all } \lambda \in \R \}.
\]
Leaves are compact if $v_2/v_1$ is rational. We are going to explore the case 
of  $v_2/v_1$ irrational, which implies that all the leaves are dense in
$T$. In this case, the leaf passing through $x$ has intrinsic 
topology $\tau_{\mathrm{int}}$ given by the real line. 
Precisely, the bijective immersion $i: \R \rightarrow 
(\mathcal{L}_{x}, \tau_{\mathrm{int}} )$ defined by $\lambda \mapsto 
(a+\lambda v_1, b+ \lambda 
v_2) 
\ \ (\bmod \ \Z^2)$ is a homeomorphism.

The family of foliated semigroups $P_t ^{(n)}$ in 
$T$ will be described via the  Stratonovich 
differential equation:

    \begin{equation}\label{Eq: do exemplo folheado no toro}
        dx_t \ = \ v\ d B_t,
    \end{equation}
where $B_t$ is the standard Brownian motion and the flow $\varphi_t (x) = 
(a + v_1 B_t, b + v_2 B_t )\bmod \Z^2$ preserves each leaf of this dense 
foliation in the sense that each trajectory lays in a single leaf. 
The family of semigroups $P_t ^{(n)}$ are given by, for 
$f \in C(T^n)$,
\[
P_t ^{(n)} f(x_1, \ldots , x_n) =  \E \left[  f(\varphi_t (x_1), 
\ldots, \varphi_t (x_n))  \right] 
\]
with $(x_1, \ldots , x_n) \in T^n$.
In particular, the first semigroup  
\[
P_t 
^{(1)} f (x) = \int_T f \left( (x+ v z) \bmod \Z^2 
\right) \ g_t(z)\ dz
\]
is foliated with semigroups in the leaves $P_t^{\mathcal{L}}$ given as follows: 
given a function $f\in C_0(\mathcal{L}_x, \tau_{\mathrm{int}})$, 
\[
 P_t^{\mathcal{L}} f(x) = \int_{\R} f \left( i (z+ i^{-1}(x)) \right)\ g_t(z) \ 
dz.
\]
Here we have used that for $z\in \R$,
\[
        g_t(z) \ = \ \frac{1}{\sqrt{2\pi t}}\exp{\left(-\frac{z^2}{2t}\right)}
    \] 
is the one dimensional heat kernel associated to linear Brownian motion, 
with $t > 0$.
In this particular example the Markovian semigroups in 
the leaves $P_t^{\bar{\mathcal{L}}}$ are Feller in $C_0 
(\bar{\mathcal{L}})$, see Remark after Corollary \ref{Cor: sg Feller na folha}. 
The 
main differences between $P_t^{}$ and 
$P_t^{\mathcal{L}}$ 
which are relevant for our technique here are: 
\begin{enumerate}
 \item The domain of $P_t^{}$ restricts to 
$\Z^2$-periodic function on $\R^2$ while the domain of 
$P_t^{\mathcal{L}}$, with $\mathcal{L} \in \mathcal{F}$ 
extends to the non-compact topology of $\R$. Intersection of the domains is the 
unitary set of the null function.
  \item  The support of the transition probability measures $\mu_{P_t}(x, dy) $ 
associated to 
$P_t$ is the whole manifold $T$; while the 
support of the probability measures $\mu_{P^\mathcal{L}_t}(x, dy) $ is the 
leaf 
$\mathcal{L}_x$.
\end{enumerate}
The second item above is precisely the phenomenon that 
the supports of the measures associated to the semigroups on the leaves 
restrict to the leaves themselves, even if the leaves are dense in $M$. This 
property has been exploited in the 
proof of Theorem \ref{Thm: fluxo_foleado}. 

\eop

We finish this section with a remark on the existence of intrinsic LJR flow on 
each leaf.

\begin{remark} \emph{ Let $( P_{t}^{(n)}: n \in \N )$ be, as before, a 
CDP-family of
Feller semigroups in $M$. Denote by $\varphi_{s,t}$ the corresponding LJR-flow 
in $M$. Introduce the foliated structure $\mathcal{F}$ in $M$ and assume 
that the family of semigroups is not only foliated but also that the semigroups 
$P_t^{\mathcal{L}}$ in the leaves $\mathcal{L} \in \mathcal{F}$ are also Feller 
(either by topological 
reasons, e.g. compact leaves, or by more general condition as in Corollary 
\ref{Cor: sg Feller na folha}). Note that diagonal preserving condition  
(2) in Definition \ref{Def: familia CDP} is trivially satisfied when beforehand 
we have a flow associated to the semigroup. Hence, in each leaf 
$\mathcal{L}$, it exists a CDP-family of Feller semigroups in $\mathcal{L}$ 
which satisfies again the hypothesis of Theorem \ref{Thm_fluxo_LJR}. 
It means that each leaf $\mathcal{L}$ has intrinsically their own LJR-flow 
$\psi_{s,t}$ defined in the probability space $\Omega^{\mathcal{L}}= 
(\Pi_{s\leq t} F)$ , where $F$ is the space of measurable mappings on 
$\mathcal{L}$ with the appropriate $\sigma$-algebra, as described in Section 
1.1.
 An alternative and natural choice of $\psi_{s,t}$ is 
$\varphi_{s,t}|_{\mathcal{L}}$, $\omega$-wise based on the same previous 
probability space $\Omega$. But in general the relation between 
$\psi_{s,t}$ and $\varphi_{s,t}$ weakens to 
the 
average: for $x \in \mathcal{L}$ and $f \in C(M)$:
\[
 \E^{\Omega^{\mathcal{L}}} \left( f|_{\mathcal{L}} \circ \psi_{s,t} (x)
        \right)  =  \E \left( f\circ
        \varphi_{s,t}(x) \right).
\]
}
\end{remark}

\subsection{Coalescing foliated flows}

In this section, we consider coalescent foliated semigroups. As before, let 
$(P_t^{(n)})$ be a foliated family of Feller semigroups in a compact 
Riemannian manifold $M$  endowed with a regular foliation $\mathcal{F}$. 

We denote by $X_t^{(n)}= (X^1_t, \ldots , X^n_t), n\geq 1$, with $X_0 ^{(n)}= 
(x_1, \ldots, x_n)$ 
the Markovian processes in $M^n$ associated to the laws of $P_t^{(n)}$ 
starting at 
the point $(x_1, \ldots, x_n)$. Consider the partial diagonals $\Delta_n = \{x 
\in M^n: \mbox{ there exists a pair } i \neq j \mbox{ with } x_i=x_j  \}$ and 
the entry times $T_{\Delta_n} = \inf \{ t\geq 0, X_t^{(n)} \in \Delta_n\} $. 
With the same notation as in \cite{Le Jan-Raimond}:

\begin{theorem}[Le Jan-Raimond] \label{Thm: coalescencia LJR} There exists a 
unique 
compatible family  of Markovian semigroups $ \{ P_t^{(n), c}, n\geq 1 \}$ on $M$
 such that if $X^{(n),\, c}$ is the associated $n$-point motion and 
$T^c_{\Delta_n} = \inf \{t \geq 0, X^{(n),c} \in \Delta_n\} $, then:
\begin{description}
 \item{(a) } $(X_t^{(n), c}, t\leq T^c_{\Delta_n})$ is equal in law to  
$(X_t^{(n) }, t\leq T_{\Delta_n})$;
\item{(b) } for $t\geq T^c_{\Delta_n}$, we have that $X_t^{(n), c} \in \Delta_ 
n$.
\end{description}

\end{theorem}

With further condition it is possible to guarantee that the coalescent 
semigroup $P_t^{(n), c}$ is Feller. In fact, denoting by 
$P_{(x,y)}^2$ the transition probability associated to $P^{(2)}_ t$ at the 
point $(x,y)$, introduce:

\bigskip

\noindent {\bf Condition C:} For all $t>0, \epsilon >0$ and $x\in M$, assume 
that 
$$
 \lim_{y \rightarrow x} P_{(x,y)}^2  \left[ \{t < T_{\Delta_2}  
\} \cap \{ d(X^1_t, X^2_t) > \epsilon \} \right] =0, 
$$
and for some pair $x,y \in M$, $P_{(x,y)}^2 [T_ {\Delta_2} < \infty] >0$.

\bigskip 

We can construct now a foliated coalescing flow in $M$:

\begin{proposition} \label{Prop: fluxo_foleado e coalescente} Let $(P_t^{(n)}, 
n\geq 1)$ be a foliated family of Feller semigroups on $(M, \mathcal{F})$ which 
satisfies Condition C above. 
Then the coalescing semigroups $P_t^{(n), c}$ are  CDP Feller 
foliated semigroups for all $n\geq 1$, hence associated to a coalescing 
foliated flow 
$\varphi_{s,t}$. 
 
\end{proposition}

\proof For every $n \geq 1$, by concatenating a Markov process which stops 
when it hits a partial diagonal set $\Delta_n$ with a process starting from 
this corresponding point, Le Jan and Raimond \cite[Thm 4.1]{Le Jan-Raimond} 
have constructed explicitly a Markov process $X^{(n), c}$ with the properties 
established by Theorem \ref{Thm: coalescencia LJR}. Additional Condition C 
implies 
also that the family $P_t^{(n), c}$ is a CDP Feller semigroup, \cite[Thm 
4.1]{Le Jan-Raimond}.

We only have to note that, by construction, for $n=1$, the law of $X^{(1), 
c}$ and $X^{(1)}$ are equal, hence $P_t^{(1)} = P_t^{(1), c}$ for all $t\geq 
0$.
This implies that $P_t^{(1), c}$ is foliated.  Hence, the result follows by 
Theorem \ref{Thm: 
fluxo_foleado} and Corollary \ref{Cor: equivalencias de folheacoes}.

\eop

Proposition \ref{Prop: fluxo_foleado e coalescente} in particular implies that 
the pair of points $x,y \in M$ such that  $P_{(x,y)}^2 [T_ {\Delta_2} < \infty] 
>0$, in Condition C, must be in the same leaf. In fact, $\varphi_{s,t}(x) \in 
\mathcal{L}_x$ and $\varphi_{s,t}(y) \in 
\mathcal{L}_y$ and the intersection $\mathcal{L}_x \cap \mathcal{L}_y$ is non 
empty if and only if $x$ and $y$ are in the same leave.

\section{An averaging principle for foliated LJR flow}

 In this section we apply the technique of  foliated semigroup and foliated flow
to obtained an averaging  principle for foliated LJR flows when the
leaves of the foliation are compact. Initially we introduce appropriate
foliated coordinates such that the leaves are going to be mapped in
horizontal plaques of Euclidean space and we will have also a coordinate system
for a direction which is transversal to the leaves. For convenience, in this
coordinate system the leaves will be called
\textit{horizontal} and the transversal direction will be called the
\textit{vertical} direction.

Given a family of CDP foliated semigroups $(P^{(n)}_t, n \in \N)$ in $(M, 
\mathcal{F})$ we 
are going to consider
a small first order perturbation in the associated LJR-flow, corresponding to 
a family of CDP semigroups $(P^{(n), \epsilon}_t, n \in \N)$, generically no 
longer 
foliated. Precisely, if $P_t^{\epsilon}$ denotes the first perturbed semigroup, 
then there exists a vector field $K$ in $M$
(generically transversal to the leaves) such that for a function $f \in
C^{\infty}(M)$,
\[
 \lim_{t \searrow 0} \frac{P^{\epsilon}_t f - P_t f}{t} = \epsilon K 
f.
\]
We also localize our hypothesis such that $P^{\epsilon}_t f (x)$ is determined 
by $f$ restricted to a
neighbourhood of the leaf $\fL_x$ for small $t$. Equivalently, the support of 
the
probability measure associated to the Feller semigroup $P^{\epsilon}_t f
(x)$ stays in a neighbourhood of $\mathcal{L}_x$ for sufficiently small $t\geq 
0$.

Our main result in this section establishes an averaging principle for the
dynamics induced by the family $P^{(n),\epsilon}$ i.e. by its LJR flow in the 
transversal component. More precisely, as
$\epsilon$ goes to zero, the average of the vertical component of the perturbed
flow $\varphi^{\epsilon}_{0,t}$ approaches in a certain topology the solution of
the ODE in the vertical space given by the average of the perturbing vector
field $K$ in each leaf,
where the average is taken according to the invariant measure in each leaf of
the unperturbed system generated by the family $P^{(n)}_t$.

The approach here generalizes to semigroups the results for continuous
diffusions \cite{ivan}, foliated L\'evy processes \cite{Hoegele-Ruffino} and 
for Hamiltonian/symplectic
structures in X.-M.-Li \cite{Li}.

\paragraph{The foliated coordinate system.}

Given an initial condition $x_0\in M$, let $U\subset M$ be a
bounded neighborhood of
$x_0$ which is diffeomorphic to $\fL_{x_0}\times V$, with $V$ a connected open 
set in $\R^d$ containing the origin. By compactness of $\fL_{x_0}$, there 
exists a finite number of local
foliated
coordinate systems $\psi_i: U_i \rightarrow W_i \times V \subset
\R^n \times \R^d $, where $W_i$ and $V$ are open sets, say with $1 \leq i \leq
k$ and $x_0 \in U_1$ such that:

\begin{description}
 \item [1)] $U = \cup_{i=1}^k U_i $;
\item [2)] The leaf $\fL_{x_0}= \cup_{i=1}^k \psi^{-1}(W_i
\times \{0\})$, i.e. each $U_i$ is
diffeomorphic to the product of an open set  in the
leaf $\fL_{x_0}$ and the vertical component  $V$;
\item [3)] If a pair of points
$p\in U_i$ and $q \in U_j$ in $U$ belong to the same leaf then their
transversal
coordinates in $V$ are the same;  i.e. $\pi (\psi_i(p)) = \pi (\psi_j(p))$
where $\pi$ is the projection on the
transversal space $V$;
\end{description}

Note that for a fixed $y\in V$, the finite
union $\displaystyle \cup_{i=1}^k \psi_i^{-1} (W_i \times \{ y\} )$ is the leaf
$\fL_{\psi_i^{-1} (x,y)}$ for any $x\in W_i $. Natural examples
of this scheme of coordinates systems appear if we consider compact foliation
given by the inverse image of
submersions: values in the image space provide
local coordinates for the vertical space $V$.

Item (3) above also allows to simplify the notation in such a way that we can 
omit
the coordinate system when dealing with the vertical directions, i.e. we shall
write $\pi: M \rightarrow V$ to denote $\pi \circ \psi_i$, independently of
the (finitely many) index $i$. In coordinates, we write 
\[
\pi(\cdot)=(\pi_1(\cdot), \ldots, \pi_d(\cdot)) \in V \subset \R^d.
\]

\bigskip

\noindent \textit{Hypotheses on the perturbed semigroup.} We are going to assume
that following behaviour in the transversal dynamics. We shall denote by
$\varphi_{s,t}$ and
$\varphi^{\epsilon}_{s,t}$ the LJR flows associated to $P_t$ and
$P^{\epsilon}_t$ respectively.
\begin{description}
 \item [H1)] Vertical regularity a.s. of the perturbed flow: For all $i=1,
\ldots , d$,
 \[
 \left. \frac{d}{dt}  \pi_i (\varphi^{\epsilon}_{0,t} (x_0)) \right|_{t=0}=
\epsilon\ d \pi_i (K) (x_0).
 \]

 \item[H2)] Transversal weak boundedness of the perturbation: Denote by
$y_{t}= \varphi_{0,t}(y_0)$ and $y^{\epsilon}_{t}=
\varphi^{\epsilon}_{0,t}(y_0)$ the trajectories of the perturbed and unperturbed
systems respectively, both starting at $y_0$. Suppose there
exists a common probability space where the random flows
$\varphi_{s,t}$ and $\varphi^{\epsilon}_{s,t}$ are based such that, for a $p\in
[1, \infty )$ and any $g\in C(M)$, there exists a positive function $h(\epsilon,
t)\geq 0$, defined for $\epsilon, t \geq 0$ which is
continuous, $h(0,t)=h(\epsilon,0)= 0$ and satisfies
 \begin{equation} \label{Hypothesis 2}
  \left[ \E \Big( \sup_{0 \leq s < t} \ | g (y_{s}) - g(y_s^{\epsilon}) 
|^p\Big)
 \right]^{\frac{1}{p}} \leq h(\epsilon, \sqrt{\epsilon} t).
 \end{equation}
\end{description}
\bigskip

Note that Hypothesis (H1) states for each trajectory a property which 
always holds in
the
average, in fact: Denoting by $A$ the infinitesimal generator of $P_t$ and by 
$A^{\epsilon}$ the infinitesimal generator of $P_t^{\epsilon}$ just note that 
the projections into the vertical coordinates $\pi_i$, $i=1, \ldots, d$, are in 
the kernel of $A$ and  $A^\epsilon \pi_i = \epsilon \,
d
\pi_i (K)$. This hypothesis is canonically satisfied by semigroups generated
by foliated stochastic differential equations with an $\epsilon$ perturbation of
the drift in the direction $K$. This can be easily verified by the fact that the
kernel of the derivative $d\pi_i $ includes the tangent spaces to
the leaves, hence, the differentiability follows by It\^o formula, cf. \cite[p.
15]{ivan}, also \cite{Hoegele-Ruffino}.

Without lost of generality we can assume that function
$h(\epsilon, t)$ in Hypothesis 
(H2) is increasing in $t$ for a fixed $\epsilon$. In fact, given such a 
function $h$ we have that $ \sup_{0\leq s
\leq t} h(\epsilon, s)$ also satisfies Hypothesis (H2). This hypothesis holds
if the semigroups $P_t$ and $P^{\epsilon}_t$ are generated by
perturbation of foliated stochastic (L\'evy) differential equations, 
\cite[Lemma 2.1]{ivan}, \cite[Prop. 2.1]{Hoegele-Ruffino} for $p\geq 2$, also  
completely
integrable stochastic Hamiltonian system, Li \cite{Li}. See also Remark 1
after Lemma \ref{lema4} for more generality.  
Another class of
examples includes perturbing vector fields
$K$ which commute with the infinitesimal generator $A$ of  $P_t$, which in this
case makes $P^{\epsilon}_t= P_t \circ K_t $, where $K_t$ is the flow of local
diffeomorphisms associated to a vector field $K$ in this class. In this case
function
\[
h(\epsilon,t)= \sqrt{\epsilon} \ t\ \sup_{x\in U} |K(x)|
\]
satisfies the inequality (\ref{Hypothesis 2}).

\subsection{Averaging functions on the leaves}

By compactness, the leaf $\fL_p$ passing through a point $p\in M$ contains the 
support of an
invariant measure $\mu_p$ for
the unperturbed semigroup $P_t$. We assume that $\mu_p$ is
ergodic. Consider a continuous function $g:M\rightarrow \mathbf{R}$. We
shall work
with the $\mu_p$-average of $g$, $Q^g: V\subset \R^d
\rightarrow \R$ defined for each leaf, i.e. if $v$ is the 
vertical coordinate of $p$,
$\pi(p)=v \in V$,  then:
\[
 Q^g (v) = \int_{\fL_p} g (x) \, d\mu_p(x).
\]
 We  assume the following hypothesis
on the invariant measures on the leaves:

\begin{description}
 \item[H3.] \textit{Regularity of} $Q^g$: For any Lipschitz continuous function
$g$
on
$M$, its corresponding average function $Q^g$
on the transversal space $V$  is Lipschitz.
\end{description}
 Hypothesis (H3) means that the
invariant measures $\mu_p$ for the
unperturbed foliated system has some weakly continuity with respect to the 
vertical component of $p$; say, for instance, locally
there is no sort of bifurcation of the horizontal foliated dynamics performed 
by $(P^{(n)}_t)$
when one varies the vertical parameter, as in  
\cite{ivan}.

\bigskip

We use the
derivative of each component of $
\pi(\cdot)=(\pi_1(\cdot), \ldots, \pi_d(\cdot)) \in V \subset \R^d.
$
to get the averages $Q^{ {d\pi_i (K)}}(x)$ of the real functions $g= d\pi_i
(K)$, $i= 1, \ldots, d$ on each leaf $\fL_x$. The proposition below gives an 
ergodic estimation of the error which occurs when one considers the average 
$Q^g$ instead of the original function $g$ in a time integration.

%
%

\begin{proposition} \label{Prop: Equiv 3.1-averaging}
 For $i= 1, 2, \ldots , d$,   $t\geq 0$ and $\epsilon >0$ let
    \[
        \delta_i (\epsilon, t) =  \int_0^t d \pi_i (K) (y^
        {\epsilon}_{\frac{r}{\epsilon}})
        -  Q^{d\pi_i(K)} \pi (y^ {\epsilon}_{\frac{r}{\epsilon}})\ dr.
           \]
We have the following estimates for the difference $\delta_i (\epsilon, 
t)$ 
\[
\bigg( \E | \delta_i (\epsilon, t) |^p \bigg)^{\frac{1}{p}} \leq \sqrt{t}
H(\epsilon, t)
\]
where $H(\epsilon, t)$ is continuous in $\epsilon,t \geq 0$ and $H(0,0)=0$.
\end{proposition}

\bigskip

\noindent \textbf{Remark:} Precisely, in terms of function $h(\epsilon, t)$ in 
Hypothesis (H2) we have that 
\[
  H(\epsilon,t)= \min \left\{ h(\epsilon,t)\sqrt{t}, C_1 \epsilon^{\frac{1}{4}},
C_2 \sqrt{\epsilon} t^{\frac{3}{2}}, C_3 \sqrt{\epsilon\ t} \right\}.
\]
for some positive constants $C_1, C_2$ and $C_3$. 

\bigskip

\proof
The proof consists of changing variables to get an integration in
the interval $[0, t/\epsilon]$ such that considering a convenient
partition of this interval we estimate by comparing in each subinterval the
average
of the flow of the original system (on the corresponding leaf) with the
average
of the perturbed flow (transversal to the leaves) using Hypothesis (H2).
For sufficiently small $\epsilon$,
we take the following assignment of
increments of our partition:
\[
\Delta t = \frac{t}{\sqrt{\epsilon}}.
\]
We consider the partition  $t_n=  n \Delta t$, for $0\leq n \leq N$
 such that
\[
0=t_0 < t_1 < \cdots < t_{N} \leq
\frac{t}{\epsilon},
\]
with $N=N(\epsilon)=[\epsilon^{-1/2}]$
where here $[x]$ denotes the integer part of $x$.

\bigskip

To simplify the notation, denote by $g(x)$ the function
$d \pi_i (K)(x)$. Hence, the first integrand can be written as the sum:

\[
\epsilon
\int^{\frac{t}{\epsilon}}_0  g(y^{\epsilon}_r)dr = \epsilon
\sum^{N-1}_{n=0}\int^{t_{n+1}}_{t_n}{g(y^{\epsilon}_r)dr}+ \epsilon
\int^{\frac{t}{\epsilon}}_{t_N} g(y^{\epsilon}_r)dr .
\]

Let $\varphi_{s,t}$ denote the LJR-flow presented in the previous section for 
the
unperturbed foliated semigroups $(P_t^{(n)})$, i.e. such that each trajectory 
stays in
a single leaf of the foliation. By triangular inequality, we divide our
calculation into four parts:

\begin{equation} \label{split}
 |\delta (\epsilon, t)|\leq |A_1| + |A_2| + |A_3|+ |A_4|,
\end{equation}
 where
\begin{eqnarray*}
 A_1 & = & \epsilon
\sum^{N-1}_{n=0}\int^{t_{n+1}}_{t_n}\left[g(y^{\epsilon}_r)-g(\varphi_{t_n,
r}(y^{
\epsilon}_{t_n}))\right]\, dr,\\
 & & \\
A_2 & = & \epsilon
\sum^{N-1}_{n=0}\left[ \int^{t_{n+1}}_{t_n} g(\varphi_{t_n,
r}(y^{
\epsilon}_{t_n})) \, dr -  \Delta t
Q^g (\pi (y^{\epsilon}_{ t_n})) \right],\\
 & & \\
A_3 & = &    \sum^{N-1}_{n=0} \epsilon \Delta t
Q^g (\pi (y^{\epsilon}_{ t_n})) - \int^{t}_{0}
Q^g(\pi(y^{\epsilon}_{\frac{r}{\epsilon}}))\ dr, \\
& & \\
 A_4 &  = & \epsilon
\int^{\frac{t}{\epsilon}}_{t_N} g(y^{\epsilon}_r)dr.
\end{eqnarray*}

We proceed by showing that each of the processes $A_1, A_2, A_3$ and $A_4$ above
tends to zero on compact intervals.

%
%

\begin{lemma}\label{lema1} Process $A_1$ converges to zero on compact intervals
when
$\epsilon$ goes to zero. More precisely, we have the following estimates on
the rate of convergence:
\[
\bigg( \mathbf{E} \left|A_1 \right|^p\bigg)^{\frac{1}{p}}\leq
 h(\epsilon, t)\, t,\]
where $h(\epsilon, t)$ is given by Hypothesis (H2).
\end{lemma}

\proof If  $ \frac{1}{p}+\frac{1}{q}=1$, by H\"{o}lder and triangular
inequalities we have that

\begin{eqnarray*}
\left(\mathbf{E}  \left|A_1 \right|^p\right)^{\frac{1}{p}}
& \leq & \epsilon
\sum^{N-1}_{n=0}
 \left(\mathbf{E} \left[
\int^{t_{n+1}}_{t_n}   \left| g(y^{\epsilon}_r)-g(\varphi_{t_n,
r}(y^{
\epsilon}_{t_n})) \right| d\,r
\right]^p\right)^{\frac{1}{p}}\\
& & \\
& \leq & \epsilon \sum^{N-1}_{n=0}  \displaystyle \left( \mathbf{E}
 \left[ \left( \int^{t_{n+1}}_{t_n}
dr\right)^{\frac{1}{q}} \left(
\displaystyle \int^{t_{n+1}}_{t_n}
\left|g(y^{\epsilon}_r)-g(\varphi_{t_n,
r}(y^{
\epsilon}_{t_n}))\right|^pdr \right)^{\frac{1}{p}}\right]^p
\right)^{\frac{1}{p}} \\
& & \\
 & \leq & \epsilon (\Delta t)^{\frac{1}{q}}\sum^{N-1}_{n=0}
\displaystyle \left( \mathbf{E}
 \bigg[ \displaystyle \Delta t \
\sup_{t_n \leq r< t_{n+1}} \left|g(y^{\epsilon}_r)-g(\varphi_{t_n,
r}(y^{
\epsilon}_{t_n}))\right|^p \bigg]
\right)^{\frac{1}{p}}\\
&&\\
& \leq& \epsilon \, \Delta t\
\sum^{N-1}_{n=0} \displaystyle \left( \mathbf{E}
 \bigg[ \displaystyle \sup_{t_n \leq r< t_{n+1}}
\left|g(y^{\epsilon}_r)-g(\varphi_{t_n,
r}(y^{
\epsilon}_{t_n})) \right|^p \bigg]
\right)^{\frac{1}{p}}
\end{eqnarray*}

Hypothesis (H2) together with the fact that the law of the flow $\varphi_{s,t}$
depends only on the difference $t-s$, imply that for each $0\leq n \leq N-1$
above, the function $g$ evaluated along trajectories of the perturbed system
compared with $g$ evaluated along the unperturbed trajectories, both starting at
$y_{t_n}^{\epsilon}$ satisfies:
\[
\bigg[ \mathbf{E}  \sup_{t_n \leq r< t_{n+1}}\  \left| g(y^{\epsilon}_r) -
g\left( \varphi_{t_n,
r}(y^{
\epsilon}_{t_n}) \right)  \right|^p \bigg]^{\frac{1}{p}} \leq  h(\epsilon,
\sqrt{\epsilon} \Delta t).
\]
Hence
\begin{eqnarray*}
 \bigg[\mathbf{E}  |A_1|^p \bigg]^{\frac{1}{p}}& \leq &  \epsilon\ \Delta t
\ N \ h(\epsilon, \sqrt{\epsilon} \Delta t)
\\
&=&   h (\epsilon, t)\ t.
\end{eqnarray*}

\eop

%
%

\begin{lemma}\label{lema2} Process $A_2$ in equation
(\ref{split}) goes to zero with the following rate of convergence:
\[
 \bigg[\mathbf{E}  |A_2|^p \bigg]^{\frac{1}{p}} \leq C_1 \ \sqrt{t}\
\epsilon^{\frac{1}{4}}.
\]
for a positive constant $C_1$.
\end{lemma}

\proof By Minkowsky inequality we have that
\begin{eqnarray*}
 \bigg[\mathbf{E}  |A_2|^p\bigg]^{\frac{1}{p}}& \leq & \epsilon\
 \left[ \mathbf{E}  \left|
\sum^{N-1}_{n=0} \left[
\int^{t_{n+1}}_{t_n}g\left(\varphi_{t_n,
r}(y^{
\epsilon}_{t_n})\right)\ dr- \Delta t
Q^g (\pi (y^{\epsilon}_{ t_n})) \right]  \right|^p\right]^{\frac{1}{p}} \\
&& \\
& \leq & \epsilon\  \sum^{N-1}_{n=0}
 \left[\mathbf{E}  \left|
\int^{t_{n+1}}_{t_n} g\left(\varphi_{t_n,
r}(y^{
\epsilon}_{t_n})\right)\ dr - \Delta t
Q^g (\pi (y^{\epsilon}_{ t_n}))   \right|^p\right]^{\frac{1}{p}} \\
&& \\
& = & \epsilon\ \Delta t  \sum^{N-1}_{n=0}
 \left[\mathbf{E}  \left|
\frac{1}{\Delta t}
\int^{t_{n+1}}_{t_n}g\left(\varphi_{t_n,
r}(y^{
\epsilon}_{t_n})\right) \ dr -
Q^g (\pi (y^{\epsilon}_{ t_n}))   \right|^p\right]^{\frac{1}{p}}. \\
\end{eqnarray*}
For all $n=0, \ldots, N-1$, the ergodic theorem implies that  
the two terms inside the
modulus converges to each other when $\Delta t$ goes to infinity.
Moreover, as in \cite[Lemma 3.2]{Li} by Markovian property and central limit
theorem,
the rate of
convergence has order $\frac{1}{\sqrt{\Delta t}}$ when $\Delta t$ goes to
infinity. Hence, for small $\epsilon$
we have
\begin{eqnarray*}
 \bigg[\mathbf{E}  |A_2|^p\bigg]^{\frac{1}{p}}& \leq & C_1 \epsilon\ N (\Delta
t)
\frac{1}{\sqrt{\Delta t}} \\
& =& C_1 \epsilon \left[ \epsilon^{-\frac{1}{2}} \right] \  \sqrt{t}
\epsilon^{-\frac{1}{4}}
\\
& \leq & C_1\  \sqrt{t} \ \epsilon^{\frac{1}{4}}.
\end{eqnarray*}

\eop

%
%

\begin{lemma}\label{lema3}  $A_3$
converges to zero when $t$ or $\epsilon$ go to $0$. We have the
following rate of
convergence:
\[
\bigg(\mathbf{E} \left|A_3 \right|^p\bigg)^{\frac{1}{p}}\leq C_2 \sqrt{\epsilon}
t^2,
\]
for a positive constant $C_2$.
\end{lemma}

\proof Consider the partition $\epsilon t_n$ of the interval $[0,t]$,
whose mesh goes to zero. Then, the sum in the expression of $A_3$ is
the Riemman sum of the integral which appears in second term. Hence, the
convergence to zero corresponds to the existence of the Riemann integral,
which is guaranteed by continuity of $\pi(y^{\epsilon}_r)$ (Hypothesis H1).


We calculate now an estimate for the rate of convergence to zero.
Let $C$ be the
Lipschitz
constant of $Q^g$. Then

\begin{eqnarray}
|A_3| &\leq & \epsilon \sum^{N-1}_{n=0}  \Delta t \sup_{\epsilon t_n < s \leq
\epsilon t_{n+1}}|Q^g(\pi (y^{\epsilon}_{ t_n}))- Q^g(\pi
(y^{\epsilon}_{\frac{ s}{\epsilon}})) |  \nonumber \\
&\leq & \epsilon (\Delta t) \ C \sum^{N-1}_{n=0}  \ 
\sup_{\epsilon t_n < s \leq
\epsilon t_{n+1}} | \pi (y^{\epsilon}_{\frac{\epsilon t_n}{\epsilon} })- 
\pi (y^{\epsilon}_{ \frac{s}{\epsilon} } ) |. \label{Ineq: proof lemma 3}
\end{eqnarray}
By Hypothesis (1) we have the following inequality which is independent of
$\epsilon$: 
\[
 |\pi (y^{\epsilon}_{\frac{u}{\epsilon} }) - \pi
(y^{\epsilon}_{\frac{v}{\epsilon} }) | \leq \sup_{x\in U} K(x) \ |u-v| 
\]
for all $u,v\geq 0$. Hence, continuing the estimates for $|A_3|$, Inequality
(\ref{Ineq: proof lemma 3}) above implies that 

\begin{eqnarray*}
|A_3 |& \leq  & C_2 (\epsilon \Delta t)^2 N  \\
 &  = & C_ 2 \left( \epsilon \frac{t}{\sqrt{\epsilon}} \right)^2  
\epsilon^{-\frac{1}{2}} \\
&=& C_2 \sqrt{\epsilon} t^2,
\end{eqnarray*}
for a positive constant $C_2$.

\eop

%
%

\begin{lemma}\label{lema4} Process $A_4$ converges to zero
with
\[
\bigg(\mathbf{E} \left|A_4 \right|^p\bigg)^{\frac{1}{p}} \leq C_3 t
\sqrt{\epsilon}.
\]

\end{lemma}

\proof Denote
\[
 C_3 = \sup_{x\in U} |g(x)|.
\]
The result follows straightforward since
\[
\epsilon
\left|\displaystyle\int^{\frac{t}{\epsilon}}_{t_N}{g(y^{\epsilon}_r)dr}\right|
\leq  C_3 \epsilon \Delta t =  C_3 t \sqrt{\epsilon}.
\]

\bigskip

\noindent Now, going back to the proof of Proposition \ref{Prop: Equiv 
3.1-averaging}. Note that
each of the four estimates of Lemmas \ref{lema1}--\ref{lema4} allows a
factorization which has a
common factor $ \sqrt{t}$ times a continuous function
which goes to zero when  $(t,\epsilon) \to 0$.
Explicitly, take
\[
 H(\epsilon,t)= \min \left\{ h(\epsilon,t)\sqrt{t}, C_1 \epsilon^{\frac{1}{4}},
C_2 \sqrt{\epsilon} t^{\frac{3}{2}}, C_3 \sqrt{\epsilon\ t} \right\}.
\]
Proposition \ref{Prop: Equiv 3.1-averaging}
now follows by inequality (\ref{split}).

\eop

%
%

\noindent \textbf{Remark:} The technique we have used to prove 
Proposition \ref{Prop: Equiv 3.1-averaging} can be extended in fact 
to a larger 
class of functions $h$ in inequality (\ref{Hypothesis 2}). Let $f:\R_{\geq 0}
\rightarrow \R_{\geq 0}$ be a continuous function with $f(0)=0$,
\[
\lim_{\epsilon \searrow 0} f(\epsilon)^{-1}= + \infty\ \ \ \ \ \mbox{ and
}\ \ \ \ \ \ \lim_{\epsilon \searrow 0} \epsilon f(\epsilon)^{-1}= 0.
\]
Then
inequality (\ref{Hypothesis 2}) of Hypothesis (H2)
can be restated for $h (\epsilon, f(\epsilon)t)$. In this case, in the proof of
Proposition \ref{Prop: Equiv 3.1-averaging} one has to consider the partition 
$\Delta t =
\frac{t}{f(\epsilon)}$, $N= [f(\epsilon)^{-1}]$ and the results follows by the
same arguments. As state before, for stochastic Hamiltonian
systems, using the Liouville coordinate systems on invariant torus, one can use
$f(\epsilon)= \sqrt{\epsilon}$, see \cite{Li}. For general
stochastic equations on foliated manifolds $f(\epsilon)= |\ln \epsilon
|^{-\frac{1}{2p}}$ satisfies this extended Hypothesis (H2) hence Proposition
\ref{Prop: Equiv 3.1-averaging} also holds in this case, see \cite[Lemma 
3.1]{ivan}.

\subsection{An averaging principle}

%
%

\begin{theorem} \label{Thm: principio da media} Assume that the 
unperturbed foliated
semigroups on $M$ satisfies hypotheses (H1), (H2) and (H3) above.
 Let  $v(t)$ be the solution of the deterministic ODE in the transversal
component $V\subset \R^n$,
\begin{eqnarray}
\frac{dv}{dt}= (Q^{ {d\pi_1 (K)}}, \ldots, Q^{d\pi_d(K)})  ( v(t) )
\end{eqnarray}
with initial condition $v(0)=\pi(x_0)=0$. Let $T_0$ be the time that
$v(t)$
reaches the boundary of $V$. Then, for all  $0<t<T_0$ we have that

\[
 \left[\mathbf{E}\left( \left|\pi
\left( y^{\epsilon}_{\frac{t}{\epsilon}}
\right)-v (t)\right|^p\right)\right]^{\frac{1}{p}}
 \leq  G(\epsilon, t)
\]
where $G(\epsilon, t)\geq 0$ is continuous for nonnegative   $\epsilon$ and $t$, 
it is
decreasing in $t$ for a fixed $\epsilon$ and $G(\epsilon,0)= G(0,t)=0$.

\end{theorem}

\bigskip

\noindent \textbf{Remark:} Precisely, in terms of function $h(\epsilon, t)$ in 
Hypothesis (H2) we have that the estimates above are given by
\[
  G(\epsilon,t)= \sqrt{t} e^{Ct}\ \min \left\{ h(\epsilon,t)\sqrt{t}, C_1 
\epsilon^{\frac{1}{4}},
C_2 \sqrt{\epsilon} t^{\frac{3}{2}}, C_3 \sqrt{\epsilon\ t} \right\}.
\]
for some positive constants $C, C_1, C_2$ and $C_3$. 

\bigskip


\proof

Most of the calculations have been done in Proposition \ref{Prop: Equiv 
3.1-averaging}. We
only have to note that for
each $i=1,2,\ldots, d$, by Hypothesis (H1), Jensen's inequality and Proposition
\ref{Prop: Equiv 3.1-averaging},  we
have

\begin{eqnarray}
\nonumber \left|\pi_i
\left( y^{\epsilon}_{\frac{t}{\epsilon}}
\right)-v_i (t)\right|&\leq& \int^{t\wedge
T^{\epsilon}}_0{\left|Q^{d\pi_i(K)}(\pi^{\epsilon}(s))-Q^{d\pi_i(K)}
(v(s))\right|ds}
+|\delta_i(\epsilon,t)|\\
\nonumber &\leq&
C_i \int^t_0 \left |\pi
\left( y^{\epsilon}_{\frac{s}{\epsilon}}
\right)-v (s) \right| ds +| \delta_i(
\epsilon,t)|,
\end{eqnarray}
where each $C_i$ is the Lipschitz constant of $Q^{d\pi_i(K)}$.
Summing up the $i$'s and using Gronwall's lemma we have, for a constant $C$:
\[
\left|\pi_i
\left( y^{\epsilon}_{\frac{t}{\epsilon}}
\right)-v_i (t)\right|  \leq e^{Ct} \sum_{i=1}^n |\delta_i(
\epsilon,t)|.
\]
And the result follows by Proposition \ref{Prop: Equiv 3.1-averaging}.

\eop

\subsection{Example:}

We present a simple example to illustrate the
framework where the averaging principle for perturbed foliated
semigroups  holds.
Consider $M= \R^3
- \{ (0,0,z), z\in R \}$
with the 1-dimension horizontal circle foliation of $M$ where the leaf
passing through a point $p=(x,y,z)$ is given by the circle
\[
L_p = \left\{ (
\sqrt{x^2+ y^2}
\cos \theta, \sqrt{x^2+ y^2} \sin
\theta, z),\ \theta \in [0,2 \pi] \right\}.
\]
For a point $p_0=(x_0,y_0,z_0)$, say with $x_0 \geq 0$
consider
the local foliated coordinates in the  neighbourhood $U= \R^3
\setminus \{(x, 0, z); x\leq 0; z\in \R \}$ given by cylindrical coordinates.
Hence, the coordinate system is defined  by $\psi:U
\subset
M \rightarrow (-\pi, \pi) \times \R_{>0}\times \R $, where $u \in (-\pi, \pi)$
is angular and $v=(r,z) \in
\R_{>0}\times \R$ is such that $\psi^{-1}: (u, v)\mapsto (r \cos u, r \sin u, z)
\in M$. In this coordinate system, the transversal projections $\pi_1$ and
$\pi_2$
correspond to the radial $r$-component and the $z$-coordinate,
respectively.

Consider the semigroup $P_t$ acting in
$C_0(M)$ given by the following: denoting a point
$\varphi(x)= (\theta, r,z)$ by its coordinates and writing the entries of  
$f\in C_0(M)$ in these coordinates, 
\[
 P_t f (x) = \frac{1}{2}\  \left\{ f (\theta+t, r, z)(1+ e^{-2t}) + f(\theta + 
\pi + t, r, z)
 (1-e^{-2t}) \right\}.
\]
The sum in the angles are taken module $2\pi$. This semigroup corresponds to a
L\'evy
flow in each circular leaf diffeomorphic to $S^1$ of the foliation which has
simultaneously
two commutative behaviour: pure rotation and Poisson jumps to the antipodal.
Its infinitesimal generator in cylindrical coordinates is given by $Af
(\theta,r,z) = \frac{\partial f}{\partial \theta} + f(\theta+ \pi, r,z) -
f(\theta, r,z)$.
See e.g. Applebaum \cite{Applebaum} or Liao \cite{Liao}.
Hence it is obviously a foliated semigroup in the foliated space $M$.
The unique invariant probability measure is the normalized Lebesgue measure in
each circle, hence Hypothesis (H3) is satisfied.

\paragraph{Adding a first order perturbation.} We
investigate the
effective behaviour of a small
transversal perturbation in the semigroup, such that the original infinitesimal
generator is perturbed by $\epsilon K$. Functions $Q^{d\pi_1(K)}$ and
$Q^{d\pi_2(K)}$ are simply the integral along each circle of the radial and
vertical components of $K$, respectively. 
Theorem \ref{Thm: principio da media} says that in the average, the transversal 
behaviour of $P_t^{\epsilon}$ is approximate by $v(t)$ where $v(t)$ is solution
starting at zero of the EDO $v'(t)= (Q^{d\pi_1(K)}(v_t),
Q^{d\pi_2(K)}(v_t))$, i.e.

\[
 P^{\epsilon}_{\frac{t}{\epsilon}} \left|\pi_i
\left( \cdot
\right)-v_i (t) \right|
 \leq  G(\epsilon, t).
\]
A class of examples appears if we consider a vector field $K$ which
commutes with $A$, say, suppose that in cylindrical coordinates $K$ is given
by:
\[
 K(\theta, r,z) = (0, \lambda_0, k_3(z))
\]
where $k_3: \R \rightarrow \R $ is a smooth function with bounded
derivative. The perturbed semigroup in this case is given by
\[
 P_t^{\epsilon} f (x) = \frac{f\big( \eta_t  \big)(1+
e^{-2t}) + f \left( \eta'_t \right) (1-e^{-2t})}{2}.
\]
where
\[
 \eta(t)= \Big( \theta+t,\  r+ \epsilon \lambda_0
t, \ \xi^\epsilon (z,t) \Big);
\]
\[
 \eta'_t= \Big( \theta+\pi+t,\  r+ \epsilon \lambda_0
t, \ \xi^\epsilon (z,t) \Big);
\]
and $\xi^\epsilon (z,t)$ is the solution of the ODE in the real line
\[
 \xi^\epsilon (z,t) = z + \int_0^t k_3(\xi^\epsilon (z,s))\ ds.
\]
Vector field $K$ commutes with the infinitesimal generator $A$, hence
hypothesis (H1) and (H2) with exponent $p=1$ are satisfied.

 \bigskip

 In this case we have the radial $d\pi_1(K)= \lambda_0$ and the vertical
$d\pi_2 (k)(\theta,r,z)= k_3(z)$, hence their average with respect to
Lebesgue measure on the leaves are $Q^{d\pi_1(K)}\equiv \lambda_0$ and
$Q^{d\pi_2 (k)}= k_3(\pi_2 (x))$. Hence the transversal
components as stated in the main Theorem \ref{Thm: principio da media} are given
by
$v(t)=(r_0+ \epsilon\ t \lambda_0, \xi^\epsilon (z,t))$ for all $t\geq 0$, if
$\lambda_0 \geq 0$, and $ 0 \leq t < \frac{r_0}{\epsilon |
\lambda_0|} $ if $\lambda_0<0$. One checks easily that
\[
 P^{\epsilon}_{\frac{t}{\epsilon}} \left|\pi_i
\left( \cdot
\right)-v_i (t) \right|=0
 \leq  G(\epsilon, t).
\]
\eop

\bigskip

\bigskip

\bigskip

\noindent {\bf Acknowledgments:} This article has been written while the
authors are visiting Humboldt University, Berlin. They would
like to express their gratitude to Prof. Peter Imkeller and his research group
for the nice and friendly hospitality. Author P.H.C. has been supported by
CNPq 149688/2010-5 and 236640/2012-7
and P.R.R. has been partially supported by FAPESP 11/50151-0 and
12/03992-1.


\begin{thebibliography}{99}

\bibitem{Applebaum} D. Applebaum -- {\it L\'{e}vy processes and stochastic
calculus,} Cambridge University press, 2004.

\bibitem{L Arnold}  L. Arnold -- {\it Random dynamical systems,}
Springer-Verlag, 1998.


\bibitem{Baxendale} P. Baxendale --  T. E. Harris's contributions to recurrent 
Markov processes and 
stochastic flows. {\it Ann. Probab.} 39 (2011), no. 2, 417-428.

\bibitem{Bertoin-Le Gall} J. Bertoin and J. F. Le Gall -- Stochastic flows 
associated to coalescent processes. {\it  Probab. Theory Related 
Fields},  126 (2003), 261-288.






\bibitem{camacho} C. Camacho and A. Lins-Neto -- {\it Geometric theory of 
foliations,} Birkh\"{a}user Boston, 1985.



\bibitem{candel} A. Candel and L. Conlon -- {\it Foliations I and II}.
Graduate Studies in Mathematics, American Mathematical Society,
1999.

\bibitem{Candel_Adv_Math} A. Candel -- The harmonic measures of Lucy Garnett. 
{\it
Adv. Math.} 176 (2003), no. 2, 187?247. 

\bibitem{Catuogno-Ledesma-Ruffino} P. Catuogno, D. Ledesma and P. Ruffino -- 
Harmonic
measures in embedded foliated manifolds. Submitted. (2012) ArXiv 1208.0629.

\bibitem{Chen_Xiang} J. Chen and  K. N. Xiang -- Natural flow not in Le 
Jan-Raimond framework. {\it Stochastics and Dynamics,} 12 (2012), no. 2,
1150014.



\bibitem{Garnett} L. Garnett,  \textit{Foliation, the ergodic theorem and
Brownian motion}.
Journal of Functional Analysis 51, (1983)pp. 285-311.



\bibitem{ivan} I. I. Gonzales-Gargate and P. R. Ruffino -- An averaging 
principle for
diffusions in foliated spaces.
Submitted (2012) ArXiv 1212.1587.

\bibitem{Hoegele-Ruffino} M. H\"ogele and P. R. Ruffino -- Averaging along 
L\'evy diffusions in foliated spaces. Preprint, Mathematics Department, Potsdam 
University 2 (2013) 10. Submitted.

\bibitem{kaimainovich} V. A. Kaimanovich,\textit{Brownian motion on foliations:
Entropy, invariant measures, mixing}
Functional Analysis and Its Applications, Vol. 22, N$^0$ 4, (1988) pp. 326-328





\bibitem{Kunita 1}  H. Kunita --  Stochastic differential equations and
stochastic flows of diffeomorphisms. In {\it \'{E}cole
d'Et\'{e} de Probabilit\'{e}s de Saint-Flour XII-1982}, 
pp. 143--303. Lecture Notes in Math. 1097, Springer-Verlag, Berlin,
1984.

\bibitem{Kunita 2}  H. Kunita -- {\it Stochastic flows and stochastic
differential equations}. Cambridge University Press, 1988.



\bibitem{Le Jan-Raimond} Y. Le Jan and O. Raimond -- Flows, coalescence and noise.
{\it  Ann. Probab.,}
32 (2004) no. 2, 1247-1315.

\bibitem{Le Jan- Raimond 2} Y. Le Jan and O. Raimond -- Flows associated to 
Tanaka's SDE. ALEA Lat. Am. 
J. Probab. Math. Stat. 1 (2006), 21-34. 

\bibitem{Li} X. M. Li -- An averaging principle for a completely
integrable stochastic Hamiltonian systems. {\it Nonlinearity}, 21 (2008)
803-822.





\bibitem{Liao} M. Liao -- {\it L\'evy processes in Lie groups}. Cambridge
University Press, 2004.

\bibitem{plante}J. F. Plante, \textit{Foliations with measure preserving
holonomy}.
Annals of  Mathematics, 102 (1975), 327-361.


\bibitem{Revuz-Yor} D. Revuz and M. Yor -- {\it Continuous
martingales and Brownian motion}. Springer-Verlag, Berlin 1999.



\bibitem{Tondeur} P. Tondeur. \textit{Foliations on Riemannian manifolds}.
Universitext, Springer Verlag, Berlin-Heidelberg-New York, 1988.


\bibitem{Tsirelson} B. Tsirelson. \textit{ Nonclassical stochastic 
flows and continuous products}. Probab. Surv. 1 (2004), 173--298.

\bibitem{Walcak} P. Walczak -- {\it Dynamics of foliations, groups and
pseudogroups}. Birkh\"auser Verlag 2004.

\end{thebibliography}
\end{document}